\documentclass[a4paper,12pt]{amsart}
\usepackage[utf8]{inputenc}
\usepackage{geometry}
\usepackage[colorlinks=true, linkcolor=blue, citecolor=blue, urlcolor=blue]{hyperref}
\savegeometry{paper}
\usepackage{cite,graphicx,pstricks,appendix,dsfont,amsmath,amssymb,amsthm,mathtools,comment}
\usepackage[shortlabels]{enumitem}
\usepackage{tikz,lipsum,lmodern}
\usepackage[most]{tcolorbox}
\usepackage{xcolor,blindtext}
\numberwithin{equation}{section}
\newtheorem{theorem}{Theorem}[section]
\newtheorem{corollary}[theorem]{Corollary}
\newtheorem{definition}[theorem]{Definition}

\newtheorem{proposition}[theorem]{Proposition}

\newtheorem{lemma}[theorem]{Lemma}

\newtheorem{remark}[theorem]{Remark}

\newcommand{\E}{{\mathbb E}}

\begin{document}
\title{On the geometry of Gaussian Mean Shift Detection}

\author{Pedro Abdalla}
\thanks{Pedro Abdalla: Department of Mathematics, University of Wisconsin-Madison}

\author{Grigoris Paouris}
\thanks{Grigoris Paouris: Department of Mathematics, University of Texas A\&M}

\author{Roman Vershynin}
\thanks{Roman Vershynin: Department of Mathematics, University of California Irvine}

\date{\today}

\begin{abstract}
We consider the problem of detecting a shift in the mean of a standard Gaussian random vector, where the direction of the shift belongs to an
arbitrary subset $T$ of the Euclidean unit sphere.

We introduce a correlation functional $\kappa(T)$
and establish a general lower bound on the detection
radius in terms of this functional. We identify
regularity conditions under which this lower bound
is optimal up to an absolute constant. As an application, we recover some well-known results in the literature such as the detection
radius for $s$-sparse unit vectors.
\end{abstract}

\maketitle
\section{Introduction}
High-dimensional hypothesis testing is a fundamental problem in mathematical statistics. In particular, detecting perturbations of a multivariate Gaussian vector has played a role in fundamental tasks such as sparse principal component analysis \cite{BerthetRigollet2013,JohnstoneLu2009}, higher criticism \cite{DonohoJin2004} and detecting a hidden or planted signal \cite{MendelsonPaourisVershynin2026,WeiWainwrightGuntuboyina2019,AriasCastroCandesDurand2011,Lugosi2017}.

In this note, we consider a basic Gaussian hypothesis testing problem: {\em can we detect a shift in the mean of a multivariate Gaussian random vector?} The difficulty of this problem naturally depends on (a) the geometry of the set of possible shifts, and (b) on whether the shift is chosen adversarially or sampled from a probability distribution.

Let's formalize these problems now. 

\subsection{Worst case detection problem}
Given a set $T\subseteq S^{n-1}$ and a radius $r>0$, we consider the following hypothesis testing problem:
\begin{equation}
\label{eq:hypothesis_testing}
\begin{aligned}
H_0 &: X\sim \mathcal{N}(0,I_n),\\
H_a &: X\sim \mathcal{N}(-rt,I_n), \qquad \text{for some fixed }t\in T.
\end{aligned}
\end{equation}
We define the {\em critical radius} $r(T)$ as the smallest $r>0$ for which there exists a measurable function (test)
$L:\mathbb{R}^n\to\{0,1\}$ such that
\begin{equation}
\label{eq:typeI_II_error}
\mathbb{P} \big(L(X)=1\mid X\sim\mathcal N(0,I_n) \big)\le0.1,\qquad
\sup_{t\in T}\mathbb{P} \big(L(X)=0\mid X\sim\mathcal N(-rt,I_n)\big)\le0.1.
\end{equation}
Equivalently, $r(T)$ is the smallest $r>0$ for which there is a Borel set $A$ with
\[
\gamma_n(A^c)\le0.1,\qquad \sup_{t\in T}\gamma_n(A+rt)\le0.1.
\]
The constant $0.1$ is arbitrary and could be replaced by any other constant. The goal is to characterize $r(T)$ in terms of the geometry of the set $T$. The supremum over $t \in T$ in \eqref{eq:typeI_II_error} ensures two things: (a) the non-adversarial nature of the shift: $t$ cannot be chosen depending on $X$; and (b) the shift must be detected in the worst possible case.

\subsection{Bayesian detection problem}
We also study a Bayesian version of the same question: for a measure $\mu$ supported on $T$, we consider the problem
\begin{equation}\label{eq:hypothesis_testing_Bayesian}
\begin{aligned}
H_0 &: X\sim\mathcal N(0,I_n),\\
H_a &: X\sim\mathcal N(-rt,I_n),\qquad t\sim\mu.
\end{aligned}
\end{equation}

The Bayesian (critical) radius $r_\mu(T)$ is defined analogously. Clearly, the two quantities are related $r(T)\ge r_{\mu}(T)$ for any measure $\mu$ supported on $T$. It is not clear whether this inequality can be reversed, up to a constant factor.

\subsection{Main results}
Our main result is to characterize $r(T)$ and $r_{\mu}(T)$ in terms of $T$ when the set $T$ and the measure $\mu$ satisfy some regularity conditions.

\begin{theorem}[Informal]
Let \(T\subset S^{n-1}\) and consider a measure $\mu$ in the set $P(T)$ of all probability measures supported on $T$. Set
\[
\kappa_{\mu}(T):=
\|\langle t,t'\rangle\|_{\psi_1} \quad \text{and} \quad \kappa(T) =\inf_{\mu \in P(T)}
\kappa_{\mu}(T),
\]
where \(t,t'\) are independent samples from \(\mu\). Then 
\[
r_{\mu}(T)\gtrsim
\frac{1}{\sqrt{\kappa_\mu(T)}} \quad \text{and} \quad r(T)\gtrsim
\frac{1}{\sqrt{\kappa(T)}}  .
\]
Moreover, if \(T\) and $\mu$ satisfy suitable regularity conditions, then
\[
r_{\mu}(T)\asymp
\frac1{\sqrt{\kappa_{\mu}(T)}}\quad \text{and} \quad r(T)\asymp
\frac1{\sqrt{\kappa(T)}}.
\]
\end{theorem}
The lower bounds follow from Proposition~\ref{prop:lower_bound_bayesian},
while the corresponding upper bounds under suitable regularity
conditions are established in Theorems~\ref{thm:main_Bayesian}
and~\ref{thm:regular}, respectively.

Although computing $\kappa(T)$ for a general set $T$ may be challenging, its value can be determined for sets with sufficient structure. As an illustration, in Section~\ref{s: sparse}
we compute $\kappa(T)$ for the set of \emph{sparse unit vectors}, recovering the known phase transition in the detection radius for this set.
For arbitrary sets $T$, we also establish lower bounds on
$\kappa(T)$ in terms of covering numbers and Gaussian width
in Proposition~\ref{prop:kappa_covering}.
These bounds are sharp in certain regimes. A complete geometric characterization of $\kappa(T)$ remains open.
\medskip

\subsection{Related literature}
The hypothesis testing problems \eqref{eq:hypothesis_testing} and \eqref{eq:hypothesis_testing_Bayesian} (and similar variants) have been extensively studied in the literature in a case-by-case basis \cite{AriasCastroBubeckLugosi2012,DonohoJin2004,WeiWainwrightGuntuboyina2019,Lugosi2017,IngsterTsybakovVerzelen2010,AriasCastroCandesPlan2011} but to the best of our knowledge there is no general result for an arbitrary measurable set $T\subset S^{n-1}$. For the lower bound  Proposition~\ref{prop:lower_bound_bayesian}, our technique relies on the classical second-moment method for
mixture likelihood ratios, which has been extensively
used in minimax hypothesis testing; see, for instance,
the monograph of Ingster and Suslina
\cite{IngsterSuslina2003}.  

Another work closely related to ours is that of Mendelson, Paouris, and Vershynin~\cite{MendelsonPaourisVershynin2026}. They consider an adversarial formulation of the problem in which the perturbation $rt$ may depend on the observed Gaussian vector $X$. More precisely, the adversary observes $X$ and then chooses a direction $t=t(X)\in T$ with the goal of making the perturbed observation $X+rt$ indistinguishable from an unperturbed Gaussian vector. In contrast, in our setting the direction $t\in T$ is fixed independently of $X$, although it may be chosen in a worst-case manner. Their main result is that the Gaussian width $w(T)$ captures the order of magnitude of $r(T)$ when $T$ satisfy some symmetry conditions. Despite this result implies that for such sets we have $r(T)\lesssim w(T)$, it can be highly suboptimal in many instances. A different machinery for our ``non-adversarial case" is needed.

\section{Lower Bounds on the Critical Radius}
\label{sec:lower-bounds}
In this section, we establish a general lower bound for the Bayesian radius $r_{\mu}(T)$ following standard techniques to estimate the likelihood ratio. A lower bound for the critical radius $r(T)$ follows immediately since $r(T)\ge r_{\mu}(T)$ for any measure $\mu$ supported on $T$. 

Let $P(T)$ be the set of probability measures supported on $T$. For a measure $\mu \in P(T)$, consider the quantity
\[
\kappa_\mu(T):=\|\langle t,t'\rangle\|_{\psi_1},
\]
where $(t,t')$ is an independent pair sampled from $\mu \otimes \mu$. To simplify the proofs, we define the $\psi_1$ norm as 
\[
\|\xi\|_{\psi_1}:=\sup_{p\ge1}\left(\frac{\mathbb E|\xi|^p}{p!}\right)^{1/p}.
\]
It should be clear this definition and the classical one via moment generating function are the same up to a multiplicative constant factor.

\begin{proposition}
\label{prop:lower_bound_bayesian}
Let $T\subset S^{n-1}$ and $\mu$ be supported on $T$. Then
\[
r_{\mu}(T)\gtrsim \frac1{\sqrt{\kappa_\mu(T)}}.
\]
\end{proposition}

\begin{proof}
Let $P_0=\mathcal N(0,I_n)$ be the standard Gaussian distribution and let
\[
P_1=\mathbb E_{t\sim\mu}\mathcal N(-rt,I_n)
\]
be the Bayesian mixture. We denote the  densities of $P_0$ and $P_1$ by $f_0$ and $f_1$, respectively. Since the set $T$ is a subset of the Euclidean sphere,
\[
\Lambda(x):=\frac{f_1(x)}{f_0(x)}
=e^{-r^2/2}\mathbb E_{t\sim\mu}e^{-r\langle x,t\rangle}.
\]
For $Z\sim\mathcal N(0,I_n)$, $\mathbb E\Lambda(Z)=1$ and, for independent
$t,t'\sim\mu$,
\begin{align*}
\mathbb E\Lambda(Z)^2
=e^{-r^2}\mathbb E_{t,t'}\mathbb E_Ze^{-r\langle Z,t+t'\rangle}=e^{-r^2}\mathbb E_{t,t'}e^{\frac{r^2}{2}\|t+t'\|_2^2}=\mathbb E_{t,t'}e^{r^2\langle t,t'\rangle}.
\end{align*}
Thus
\begin{align}
\label{eq:TV_second_moment}
\int_{\mathbb R^n}|f_1-f_0|
&=\mathbb E|\Lambda(Z)-1| \nonumber \le\left(\mathbb E|\Lambda(Z)-1|^2\right)^{1/2}\nonumber=\left(\mathbb E_{t,t'}e^{r^2\langle t,t'\rangle}-1\right)^{1/2}.
\end{align}
For any $\eta\in(0,1)$,
\begin{align*}
\mathbb E\exp\left(\frac{\eta}{2K}\langle t,t'\rangle\right)\le \sum_{\ell=0}^\infty
\left(\frac{\eta}{2K}\right)^\ell
\frac{\mathbb E|\langle t,t'\rangle|^\ell}{\ell!}\le\sum_{\ell=0}^\infty(\eta/2)^\ell
=\frac1{1-\eta/2}\le1+\eta.
\end{align*}
Hence if $r^2\le\eta/(2\kappa_{\mu}(T))$, then
\[
\int_{\mathbb R^n}|f_1-f_0|\le\sqrt{\eta}.
\]
If a test $L$ satisfies \eqref{eq:typeI_II_error} then for $A=\{x:L(x)=1\}$, we have that the total variation distance between $P_0$ and $P_1$ is
\[
\frac12\int_{\mathbb R^n}|f_1-f_0|\ge P_1(A)-P_0(A)\ge 0.8.
\]
We reach a contradiction as the left-hand side is at most $\sqrt{\eta}/2<1/2$.
\end{proof}
Since $r(T)\ge r_{\mu}(T)$ for every measure $\mu$ supported on $T$, it makes sense to consider the quantity
\begin{equation}
\label{eq:def_kappa}
\kappa(T):=\inf_{\mu\in\mathcal P(T)}\kappa_\mu(T)
=\inf_{\mu\in\mathcal P(T)}\|\langle t,t'\rangle\|_{\psi_1}.
\end{equation}
Clearly, $r(T)\gtrsim 1/\sqrt{\kappa(T)} $.

\subsection{The case of sparse vectors} \label{s: sparse}
In this section, we estimate the quantity $\kappa(T)$ for the case when $T=T_s$, the set of $s$-sparse vectors.  We prove the following  phase transition
\begin{theorem}
\label{thm:phase_transition_Ts}
Let $T_s$ be the set of $s$-sparse vectors in the sphere $S^{n-1}$. Then
\begin{equation}
\label{eq:phase_transition_Ts}
\kappa(T_s) \asymp
\begin{cases}
    \frac{1}{s\log(ed/s^2)}, \quad s\ll \sqrt{n}.\\
    \frac{1}{\sqrt{n}}, \quad s\gtrsim \sqrt{n}.
\end{cases}
\end{equation}
\end{theorem}
It is well known that $r(T_s)$ undergoes the same phase transition \cite{DonohoJin2004}, thus in this case our lower bound is optimal up to an absolute constant.
\begin{lemma}[Lower Bounds]
Let $T_s$ be the set of $s$-sparse vectors in the sphere. Then
\[
\kappa(T_s)\gtrsim \frac{1}{s\log(en/s^2)},
\]
for any $s\le \sqrt{n}$. On the other hand, if $s>\sqrt{n}$ we have 
\[
\kappa(T_s)\gtrsim \frac{1}{\sqrt{n}}.
\]
\end{lemma}

\begin{proof}
Fix an arbitrary probability measure $\mu$ supported on $T_s$, and let $t,t'$ be independent random vectors distributed according to $\mu$. We prove that
\[
\|\langle t,t'\rangle\|_{\psi_1}
\gtrsim \frac{1}{s\log(en/s^2)}.
\]
Fix an integer $k\ge 1$. By independence between $t$ and $t'$, we have
\[
\mathbb{E}\langle t,t'\rangle^{2k}
=
\left\|\mathbb{E}t^{\otimes 2k}\right\|_F^2.
\]
Keeping only the diagonal entries of the tensor, we obtain that
\begin{equation}
\label{eq:sparse_moment_diagonal}
\mathbb{E}|\langle t,t'\rangle|^{2k}
\ge
\sum_{i=1}^n\left(\mathbb{E}t_i^{2k}\right)^2.
\end{equation}

Since $t$ is supported on at most $s$ coordinates and $\|t\|_2=1$, the comparison between the $\ell_{2k}$ and $\ell_2$ norms gives
\[
\sum_{i=1}^n |t_i|^{2k}
\ge
s^{1-k}\left(\sum_{i=1}^n t_i^2\right)^k
=
s^{1-k}.
\]
Taking expectation,
\[
\sum_{i=1}^n\mathbb{E}|t_i|^{2k}
\ge s^{1-k}.
\]
Consequently, by the Cauchy--Schwarz inequality,
\[
\sum_{i=1}^n\left(\mathbb{E}|t_i|^{2k}\right)^2
\ge
\frac{1}{n}
\left(
\sum_{i=1}^n\mathbb{E}|t_i|^{2k}
\right)^2
\ge
\frac{1}{n s^{2k-2}}.
\]
Combining this estimate with \eqref{eq:sparse_moment_diagonal}, we conclude that
\begin{equation}
\label{eq:sparse_high_moment}
\mathbb{E}|\langle t,t'\rangle|^{2k}
\ge
\frac{1}{n s^{2k-2}}.
\end{equation}
By the definition of the $\psi_1$ norm,
\[
\|\langle t,t'\rangle\|_{\psi_1}
\ge
\left(
\frac{\mathbb{E}|\langle t,t'\rangle|^{2k}}
{(2k)!}
\right)^{1/(2k)}.
\]
Using \eqref{eq:sparse_high_moment},
\[
\|\langle t,t'\rangle\|_{\psi_1}
\ge
\frac{n^{-1/(2k)}s^{-1+1/k}}
{((2k)!)^{1/(2k)}}.
\]
Since $((2k)!)^{1/(2k)}\lesssim k$, we obtain that
\begin{equation}
\label{eq:sparse_psi_lower_k}
\|\langle t,t'\rangle\|_{\psi_1}
\gtrsim
\frac{1}{ks}
\left(\frac{s^2}{n}\right)^{1/(2k)}.
\end{equation}
We optimize for
\[
k=\left\lceil\log\left(\frac{en}{s^2}\right)\right\rceil.
\]
Assume first that $s\le\sqrt n$, thus
\[
\left(\frac{s^2}{n}\right)^{1/(2k)}
=
\exp\left(
-\frac{\log(n/s^2)}{2k}
\right)
\gtrsim 1,
\]
while
\[
k\asymp \log\left(\frac{en}{s^2}\right).
\]
It follows from \eqref{eq:sparse_psi_lower_k} that
\[
\|\langle t,t'\rangle\|_{\psi_1}
\gtrsim
\frac{1}{s\log(en/s^2)}.
\]

Since $\mu$ was an arbitrary probability measure supported on $T_s$, taking the infimum over measures supported in $T_s$ finishes the proof for $s\le \sqrt{n}$. The case $s>\sqrt{n}$ follows a similar proof with the choice of $k=1$. 
\end{proof}

For the corresponding upper bounds in \eqref{eq:phase_transition_Ts}, recalling  \eqref{eq:def_kappa}, it remains to construct a suitable prior $\mu \in P(T_s)$ witnessing the lower bounds derived above.
\begin{proof}[Proof of Theorem~\ref{thm:phase_transition_Ts}]
We divide the construction in two regimes.    

\textbf{Sparse regime: $s\ll\sqrt{n}$.}
We prove that
\begin{equation*}
    \kappa(T_s) \lesssim \frac{1}{s\log(ed/s^2)},
\end{equation*}
which gives a sharp upper bound when $s\ll \sqrt{n}$. To this end, pick $I$ uniformly among all subsets of $[n]$ of cardinality $s$ and consider i.i.d. random signs $\varepsilon_1,\cdots,\varepsilon_n$. Set $t$ to be 
\begin{equation*}
    t=\frac{1}{\sqrt{s}}\sum_{i\in I}\varepsilon_ie_i.
\end{equation*}
Clearly, $t\in T_s$. For an independent copy $t'$ we have
\begin{equation*}
\langle t,t'\rangle = \frac{1}{s}\sum_{i\in I\cap I'}\varepsilon_i\varepsilon_i'
\end{equation*}
Let $K=|I\cap I'|$ be the cardinality of the overlap and notice that 
\begin{equation*}
    |\langle t,t'\rangle| \le \frac{K}{s}.
\end{equation*}
Conditioning on the set $I$, $K$ becomes a hyper-geometric distribution with parameters $(n,s,s)$. By Hoeffding's inequality we have that for any $\lambda>0$
\begin{equation*}
    \E e^{\lambda K}\le \E e^{\lambda K_{bin}},
\end{equation*}
where $K_{bin}$ is a Binomial random variable with parameters $(s,s/n)$. By the standard MGF estimates for binomial random variables we have that
\begin{equation*}
\E e^{\lambda K}\le \E e^{\lambda K_{bin}}\le e^{\frac{s^2}{n}(e^{\lambda}-1)}.
\end{equation*}
Consequently replacing $\lambda \rightarrow 1/\lambda s$ we have
\begin{equation*}
    \E^{|\langle t,t'\rangle|/\lambda}\le \E^{ K/\lambda s}\le \E^{\lambda K_{bin}/\lambda s}  = e^{\frac{s^2}{n}(e^{1/\lambda s}-1)}.
\end{equation*}
If there exists a positive $\lambda$ that satisfies
\begin{equation*}
    e^{\frac{s^2}{n}(e^{1/\lambda s}-1)} \le 2,
\end{equation*}
then $\|\langle t,t'\rangle\|_{\psi_1}\le \lambda $. It suffices to find a $\lambda$ for which
\begin{equation*}
    \frac{s^2}{n}(e^{1/\lambda s}-1)\le \log 2,
\end{equation*}
which translates to
\begin{equation*}
    \lambda = \frac{1}{s\log(\frac{n}{s^2}\log 2+1)}.
\end{equation*}
\textbf{Dense regime: $s\gtrsim\sqrt{n}$.} In this case, we need to prove that $\kappa(T_s)\lesssim 1/\sqrt{n}$.

\medskip

Since $P(T_{c\sqrt{n}})\subset P(T_s)$ for some constant $c>0$ we have that 
\begin{equation*}
    \kappa(T_{c\sqrt{n}}) \ge \kappa(T_{s}).
\end{equation*}
By repeating the same proof as in the sparse regime, obtain that 
\begin{equation*}
    \kappa(T_{c\sqrt{n}})\lesssim \frac{1}{s\log(\frac{n}{s^2}\log 2+1)}\bigg
    |_{s=c\sqrt{n}} \lesssim \frac{1}{\sqrt{n}}.
\end{equation*}
Combining the two last inequalities, we obtain that $\kappa(T_s)\lesssim 1/\sqrt{n}$ whenever $s\ge c\sqrt{n}$ for some constant $c>0$. 
\end{proof}
\subsection{Connections with Convex Geometry}
It is desirable to estimate $\kappa(T)$ in terms of
geometric properties of $T$. In this section, we make
progress in this direction by establishing a lower bound
on $\kappa(T)$ in terms of the covering numbers
$N(T,\varepsilon)$, defined as the minimal number of
closed balls of radius $\varepsilon$ needed to cover $T$.
As a consequence, we also obtain a lower bound in terms
of the Gaussian width $w(T)$.

\begin{proposition}
\label{prop:kappa_covering}
Let $T\subseteq S^{n-1}$ and $\kappa(T)$ as in \eqref{eq:def_kappa}. Then
\[
\kappa(T) \gtrsim
\sup_{\varepsilon\in(0,\sqrt{2})}
\frac{2-\varepsilon^2}
{1+\log N(T,\varepsilon/2)}.
\]
In particular, 
\[
\kappa(T) \gtrsim  \frac{1}{1+w^2(T)}.
\]
\end{proposition}
The estimate cannot be improved in general. For example, for $T=T_s$ with $s=n^{\beta}$ with $\beta<1/2$ we have 
\[
\kappa(T_s) \asymp \frac{1}{s\log(en/s)}\asymp \frac{1}{w^2(T)}.
\]
However, it does not fully capture the order of $\kappa(T)$ as 
\[
\kappa(S^{n-1}) \asymp \frac{1}{\sqrt{n}} \gg \frac{1}{n} = \frac{1}{w^2(S^{n-1})}.
\]
It remains an open question the complete geometric characterization of $\kappa(T)$.

\begin{proof}
By Sudakov's inequality, for any $\varepsilon>0$ we have
\begin{equation*}
    \varepsilon\sqrt{\log N(T,\varepsilon)}\lesssim  w(T).
\end{equation*}
Fix $\varepsilon\in(0,\sqrt{2})$ and a probability measure $\mu$ supported on $T$. Let
\[
N=N(T,\varepsilon/2).
\]
Since $T$ can be covered by $N$ balls of radius $\varepsilon/2$, at least one of these balls, say $B$, satisfies
\[
\mu(B)\ge \frac{1}{N}.
\]
Moreover, $B$ has diameter at most $\varepsilon$. Hence, if $t,t'\in B$, then
\[
\langle t,t'\rangle
=
1-\frac{1}{2}\|t-t'\|_2^2
\ge
1-\frac{\varepsilon^2}{2}.
\]
Consequently, by independence of $t$ and $t'$,
\[
\mathbb{P}\left(
|\langle t,t'\rangle|
\ge
1-\frac{\varepsilon^2}{2}
\right)
\ge
\mathbb{P}(t,t'\in B)
=
\mu(B)^2
\ge
\frac{1}{N(T,\varepsilon/2)^2}.
\]
Therefore, by standard characterization for the $\psi_1$ norm (see, e.g, \cite[Proposition 2.8.1]{Vershynin_2026}) there exists $C>0$ such that for any $K>C\|\langle t,t'\rangle\|_{\psi_1}$, we have
\begin{align*}
2
&\ge
\mathbb{E}\exp\left(\frac{|\langle t,t'\rangle|}{K}\right)\\
&\ge
\exp\left(\frac{1-\varepsilon^2/2}{K}\right)
\mathbb{P}\left(
|\langle t,t'\rangle|
\ge
1-\frac{\varepsilon^2}{2}
\right)\\
&\ge
\exp\left(\frac{1-\varepsilon^2/2}{K}\right)
\frac{1}{N(T,\varepsilon/2)^2}.
\end{align*}
Taking logarithms gives
\[
\frac{1-\varepsilon^2/2}{K}
\le
\log 2+2\log N(T,\varepsilon/2).
\]
Thus
\[
K
\ge
\frac{1-\varepsilon^2/2}
{\log 2+2\log N(T,\varepsilon/2)}.
\]
Taking the infimum over $K$, we have
\[
\|\langle t,t'\rangle\|_{\psi_1}
\gtrsim
\frac{1-\varepsilon^2/2}
{1+\log N(T,\varepsilon/2)}.
\]
Since this holds for every probability measure $\mu$ supported on $T$, we conclude the proof by taking supremum over $\varepsilon \in (0,\sqrt{2})$.
\end{proof}

\section{Critical Radius: The Bayesian Setting}
\label{sec:bayesian}
In this section, we prove that the lower bound
$r_{\mu}(T)\gtrsim1/\sqrt{\|\langle t,t'\rangle\|_{\psi_1}}$
is sharp whenever $\mu$ satisfies a suitable regularity condition.

The regularity condition relates two complementary properties of the prior \(\mu\). The first is a coordinate small-ball property: a typical vector sampled from \(\mu\) has at least \(m\) coefficients of magnitude at least \(u\) in a fixed orthonormal family. The second concerns the correlation between independent samples \(t,t'\sim\mu\), measured by \(\kappa_\mu(T)=\|\langle t,t'\rangle\|_{\psi_1}\). Small values of \(\kappa_\mu(T)\) indicate that independent signals are weakly correlated, making the corresponding Gaussian mixture difficult to distinguish from the null. The regularity condition requires that the coordinate magnitude \(u\) be sufficiently large relative to this correlation scale. Consequently, at the radius \(r\asymp\kappa_\mu(T)^{-1/2}\), a typical signal produces sufficiently strong shifts in at least \(m\) coordinates for a counting test to distinguish it from Gaussian noise.

\begin{definition}[Regularity]
\label{def:regular_bayes}
A measure $\mu$ supported on $T\subset S^{n-1}$ is $(m,\delta)$-regular if there
exist an orthonormal family $v_1,\ldots,v_N$, $u\in(0,1]$, and
$1\le m\le N\le n$ such that
\[
\mu\left(\left\{t:\#\{j:|\langle t,v_j\rangle|\ge u\}\ge m\right\}\right)
\ge1-\delta,
\]
and
\[
\kappa_\mu(T)\leq C\frac{u^2}{1+\log_+(N/m^2)}.
\]
Here $C>0$ is an absolute constant.
\end{definition}

Define
\[
M(X)=\sum_{j=1}^N\mathds{1}\left(|\langle X,v_j\rangle|\ge\frac{ru}{2}\right),
\qquad
p=\mathbb P(|g|\ge ru/2),
\]
and
\[
L(X)=
\begin{cases}
0,&M(X)\le Np+m/4,\\
1,&\text{otherwise}.
\end{cases}
\]
\begin{theorem}
\label{thm:main_Bayesian}
Let $\mu$ be an $(m,\delta)$ regular measure for some $\delta < 0.1$. There exists a positive constant $C_r(\delta)$ depending only on $\delta$ such that
\[
r_{\mu}(T)\le \frac{C_r}{\sqrt{\kappa_\mu(T)}}.
\]
\end{theorem}
\begin{remark}
The first condition in Definition~\ref{def:regular_bayes}
is satisfied by every probability measure $\mu$
supported on $T\subset S^{n-1}$ for some choice
of $m$, $u$, and $N$. For example, by taking $N=n$,
$m=1$, $u=n^{-1/2}$, and the canonical basis,
the condition holds with $\delta=0$.
Without the second condition in our regularity definition,  the same proof yields
\[
r_{\mu}(T)
\leq
\frac{C(\delta)}{u}
\sqrt{1+\log_+(N/m^2)},
\]
provided the first condition holds for some
orthonormal family $v_1,\ldots,v_N$,
$1\leq m\leq N$, and $u\in(0,1]$.
\end{remark}
\begin{proof}
\textbf{Under $H_0$.}
Since $v_1,\cdots,v_N$ is an orthonormal family and $X$ is a standard Gaussian, we have that $M(X)\sim\operatorname{Binom}(N,p)$. The Gaussian tail and regularity imply
\[
\frac{Np}{m^2}\le\frac{2N}{m^2}e^{-r^2u^2/8}=o_{C_r}(1).
\]
Since $\mathbb EM=Np$ and $\operatorname{Var}(M)\le Np$, Chebyshev gives
\[
\mathbb P\left(M\ge Np+\frac m4\right)
\le\frac{16Np}{m^2}=o_{C_r}(1).
\]

\medskip

\textbf{Under $H_a$.}
We first condition on $t$. Thus,
\[
\langle X,v_j\rangle=g_j+r\langle t,v_j\rangle,
\]
where the $g_j$'s are independent standard Gaussians. Set $J(t)$ to be the set of regular points,
\[
J(t)=\{j:|\langle t,v_j\rangle|\ge u\}.
\]
For any index $j\in J(t)$, notice that
\[
\mathbb P\left(|\langle X,v_j\rangle|\ge\frac{ru}{2}\mid t\right)
\ge1-\frac p2.
\]
Moreover, for every $a\in\mathbb R$,
\[
\mathbb P(|g+a|\ge ru/2)\ge\mathbb P(|g|\ge ru/2)=p.
\]
Thus for every regular $t$, $M(X)\mid t$ stochastically dominates
$B=B_1+B_2$, where $B_1$ and $B_2$ are independent and distributed as
\[
B_1\sim\operatorname{Binom}(m,1-p/2),\qquad
B_2\sim\operatorname{Binom}(N-m,p).
\]
Next, by linearity of expectation
\[
\mathbb EB=Np+m\left(1-\frac{3p}{2}\right).
\]
Consequently, by taking $C_r$ large enough that $p<1/6$ gives
\[
\mathbb EB\ge Np+\frac{3m}{4}.
\]
Also
\[
\operatorname{Var}(B)
=(N-m)p(1-p)+m(1-p/2)p/2\le\frac32Np.
\]
Therefore
\[
\mathbb P\left(M<Np+\frac m4\mid t\right)
\le\frac{6Np}{m^2}.
\]
Averaging over $t$ and using regularity,
\[
\mathbb P\left(M<Np+\frac m4\right)
\le\frac{6Np}{m^2}+\delta.
\]
For sufficiently small absolute $\delta$ and sufficiently large $C_r$, this is
at most $0.1$.
\end{proof}

\subsection{Application to Sparse Vectors}
\label{subsec:bayesian-sparse}
We apply Theorem~\ref{thm:main_Bayesian} to the set $T_s$ of $s$-sparse vectors in the unit sphere. This illustrates that our regularity condition captures the whole phase transition.

\medskip
\textbf{Sparse Regime.}
Assume $s\ll\sqrt n$. Let $\mu_s$ be uniform on the flat $s$-sparse vectors
\[
S=\left\{\frac1{\sqrt s}\sum_{i\in I}\varepsilon_i e_i:
|I|=s,\ \varepsilon_i\in\{\pm1\}\right\}.
\]
Choose $v_i=e_i$, $N=n$, $m=s$, and $u=1/\sqrt s$. Then the first part of the regularity
condition holds with $\delta=0$. Moreover, an analogous argument used to derive the estimates for $\kappa(T_s)$ implies that
\[
\kappa_{\mu_s}\asymp\frac1{s\log(en/s^2)}
=\frac{u^2}{\log(en/m^2)},
\]
so the second part of the regularity condition also holds.

\medskip
\textbf{Dense Regime.}
Assume $s\gtrsim\sqrt n$. Let $s_0=c\sqrt n$ for a sufficiently small absolute
constant $c>0$, and let $\mu_{s_0}$ be uniform on the flat $s_0$-sparse vectors.
Since $s_0\le s$, this prior is supported on $T_s$. Taking $v_i=e_i$,
$N=n$, $m=s_0$, and $u=s_0^{-1/2}$ verifies the first regularity condition with
$\delta=0$. Furthermore,
\[
\frac{u^2}{\log(en/m^2)}
=\frac1{s_0\log(en/s_0^2)}
\asymp\frac1{\sqrt n},
\]
while it is not hard to see that 
\[
\kappa_{\mu_{s_0}}(T_s)\asymp n^{-1/2}.
\]

\section{Critical Radius: Worst-Case Detection}\label{sec:nonadv}

In view of the lower bound $r(T)\gtrsim 1/\sqrt{\kappa(T)}$ of Proposition~\ref{prop:lower_bound_bayesian}, we ask under which conditions this estimate is sharp. The Bayesian analysis of Section~\ref{sec:bayesian} suggests the answer. There, the count test succeeds because $\mu$-most points of $T$ have at least $m$ coefficients of size at least $u$ with respect to an orthonormal family. In the worst-case problem, the test has to succeed against \emph{every} $t\in T$, including points that $\mu$ may not put any weight.

\medskip
Throughout this section, $g$ denotes a standard Gaussian random variable and $\varphi$ its density. We fix an orthonormal family $V=\{v_1,\dots,v_N\}\subset\mathbb{R}^n$. For $x\in\mathbb{R}^n$, we write
\[
|x|_{(1)}\ge |x|_{(2)}\ge\cdots\ge |x|_{(N)}
\]
for the non-increasing rearrangement of $|\langle x,v_1\rangle|,\dots,|\langle x,v_N\rangle|$. We fix once and for all the absolute constants
\begin{equation}\label{eq:constants}
C_1:=7,\qquad c_1:=\frac{1}{1600},\qquad K:=4+4\log 4+2\log(1/c_1).
\end{equation}
Let $a_-<a_+$ be defined by $\mathbb{P}(|g|\ge a_-)=c_1$ and $\mathbb{P}(|g|\ge a_+)=c_1/16$.

\subsection{The multiscale count test}

Let
\[
\mathcal{M}:=\{2^k:\ k\ge 0,\ 2^k\le\sqrt N\}
\]
be the set of dyadic scales. For $m\in\mathcal{M}$, let $a_m>0$ be the threshold defined by
\begin{equation}\label{eq:threshold}
p_m:=\mathbb{P}(|g|\ge a_m)=c_1\frac{m^4}{N^2}.
\end{equation}
Since $m\le \sqrt N$, we have $p_m\le c_1<1$, so $a_m$ is well defined. Define the count statistics
\[
M_m(X):=\#\big\{j\in[N]:\ |\langle X,v_j\rangle|\ge a_m\big\},\qquad m\in\mathcal{M},
\]
and the test
\begin{equation}\label{eq:count-test}
L(X)=\begin{cases}
1, & \text{if } M_m(X)\ge Np_m+\dfrac m4 \text{ for some } m\in\mathcal{M},\\[1ex]
0, & \text{otherwise.}
\end{cases}
\end{equation}

We record two elementary facts about the thresholds. First, since $\mathbb{P}(|g|\ge a)\le e^{-a^2/2}$ we have, $m\in\mathcal{M}$,
\begin{equation}\label{eq:threshold-size}
a_m^2\le 2\log\frac{1}{p_m}=4\log\frac{N}{m^2}+2\log\frac{1}{c_1}.
\end{equation}
Second, let $m^*$ be the largest element of $\mathcal{M}$. Then $\sqrt N/2<m^*\le\sqrt N$, so $c_1/16< p_{m^*}\le c_1$ and therefore
\begin{equation}\label{eq:top-scale}
a_-\le a_{m^*}< a_+.
\end{equation}
At the top scale, the threshold is an absolute constant.

\begin{lemma}[Type I error]\label{lem:null}
If $X\sim N(0,I_n)$, then $\mathbb{P}(L(X)=1)\le 0.1$.
\end{lemma}

\begin{proof}
Since $v_1,\dots,v_N$ are orthonormal, $\langle X,v_1\rangle,\dots,\langle X,v_N\rangle$ are independent standard Gaussians, and therefore $M_m(X)\sim\mathrm{Binom}(N,p_m)$. By Chebyshev's inequality and \eqref{eq:threshold},
\[
\mathbb{P}\Big(M_m(X)\ge Np_m+\frac m4\Big)\le \frac{16Np_m}{m^2}=16c_1\frac{m^2}{N}.
\]
Since $\sum_{m\in\mathcal{M}}m^2\le\frac43 (m^*)^2\le \frac43 N$, a union bound gives
$\mathbb{P}(L(X)=1)\le \frac{64}{3}c_1<0.1$.
\end{proof}

\subsection{Visible points and regular sets}

Fix $t\in S^{n-1}$, $r>0$ and $m\in\mathcal{M}$. For $j\in[N]$, let
\[
q_j=q_j(t,r,m):=\mathbb{P}\big(|g+r\langle t,v_j\rangle|\ge a_m\big)
\]
be the probability that the $j$-th coordinate of the signal $rt$, observed in Gaussian noise, exceeds the threshold $a_m$. For a set of indices $J\subseteq[N]$ we define the \emph{excess count} and the \emph{fluctuation}
\begin{equation}\label{eq:excess}
\Delta_{m,J}(t,r):=\sum_{j\in J}\big(q_j-p_m\big),\qquad
\sigma^2_{m,J}(t,r):=Np_m+\sum_{j\in J}q_j(1-q_j).
\end{equation}
In words, $\Delta_{m,J}$ is the expected number of additional threshold exceedances that the signal produces in the directions $J$, and $\sigma_{m,J}$ bounds the standard deviation of the count, including the null fluctuation $\sqrt{Np_m}$.

\begin{definition}[Visible points]\label{def:visible}
A point $t\in S^{n-1}$ is \emph{visible at radius $r$} (with respect to $V$) if there exist a scale $m\in\mathcal{M}$ and a set $J\subseteq[N]$ such that
\begin{equation}\label{eq:visible}
\Delta_{m,J}(t,r)\ \ge\ \max\Big\{\frac m2,\ C_1\,\sigma_{m,J}(t,r)\Big\}.
\end{equation}
\end{definition}

\begin{definition}[Regular sets]\label{def:regular}
Let $C_r>0$. A set $T\subseteq S^{n-1}$ is \emph{$C_r$-regular} (with respect to $V$) if every $t\in T$ is visible at radius
\[
r_T:=\frac{C_r}{\sqrt{\kappa(T)}}.
\]
\end{definition}

Visibility asks that, at some threshold, the signal raises the count of exceedance by at least $m/2$ and by more than the count fluctuates. This can happen in two ways. A few coordinates may be so large that each of them exceeds the threshold with probability close to one (Lemma~\ref{lem:strong}). Alternatively, many coordinates may each be only slightly more likely to exceed a constant threshold, so that their small contributions add up. The freedom to choose $J$ allows one to discard coordinates that do not help.

\begin{theorem}\label{thm:regular}
Let $T\subseteq S^{n-1}$ be $C_r$-regular. Then 
\[
 r(T)\ \le\ \frac{C_r}{\sqrt{\kappa(T)}}.
\]
\end{theorem}

\begin{proof}
The type I error is controlled by Lemma~\ref{lem:null}. It remains to bound the type II error.

Fix $t\in T$ and let $X\sim N(-r_Tt,I_n)$. Since $-X\sim N(r_Tt,I_n)$ and the statistics $M_m$ depend only on $|\langle X,v_j\rangle|$, we may assume $X=Z+r_Tt$ with $Z\sim N(0,I_n)$. Let $m\in\mathcal{M}$ and $J\subseteq[N]$ be as in Definition~\ref{def:visible}, and write $\Delta=\Delta_{m,J}(t,r_T)$ and $\sigma=\sigma_{m,J}(t,r_T)$.

The indicators $Y_j:=\mathbf 1(|\langle X,v_j\rangle|\ge a_m)$ are independent Bernoulli variables with parameters $q_j$. By Anderson's inequality, translating a centered Gaussian decreases the probability of the symmetric interval $(-a_m,a_m)$, so $q_j\ge p_m$ for every $j$. Hence, on a suitable probability space, we can couple each $Y_j$ with $j\notin J$ to a Bernoulli$(p_m)$ variable $Y_j'\le Y_j$, independently over $j$. Then
\[
M_m(X)\ \ge\ B:=\sum_{j\in J}Y_j+\sum_{j\notin J}Y'_j,
\]
where
\[
\mathbb{E}B=Np_m+\Delta,\qquad \operatorname{Var}(B)=\sum_{j\in J}q_j(1-q_j)+(N-|J|)p_m(1-p_m)\le\sigma^2.
\]
If the test accepts, then in particular $M_m(X)<Np_m+m/4$. Since $m/4\le\Delta/2$, this forces $B-\mathbb{E}B<-\Delta/2$. By Chebyshev's inequality and \eqref{eq:visible},
\[
\mathbb{P}(L(X)=0)\le \mathbb{P}\Big(B-\mathbb{E}B<-\frac{\Delta}{2}\Big)\le\frac{4\sigma^2}{\Delta^2}\le\frac{4}{C_1^2}=\frac{4}{49}<0.1.
\]
Since $t\in T$ was arbitrary, the proof is complete.
\end{proof}

\subsection{Strong coordinates and the connection with Section~\ref{sec:bayesian}}

We first collect some properties of the gain function
\[
G_a(b):=\mathbb{P}(|g+b|\ge a)-\mathbb{P}(|g|\ge a),\qquad a>0,\ b\in\mathbb{R},
\]
so that $q_j-p_m=G_{a_m}(r\langle t,v_j\rangle)$.

\begin{lemma}\label{lem:gain}
Let $a>0$.
\begin{enumerate}
\item[(i)] $G_a$ is even, $G_a(0)=0$, and $G_a(b)$ is strictly increasing in $|b|$. In particular, $G_a\ge 0$.
\item[(ii)] If $|b|\ge 2a$, then $\mathbb{P}(|g+b|\ge a)\ge 1-\mathbb{P}(|g|\ge a)$.
\item[(iii)] There exists an absolute constant $\gamma_0>0$ such that
\[
G_a(b)\ge \gamma_0\min\{b^2,1\}\qquad\text{for all } a\in[a_-,a_+],\ b\in\mathbb{R}.
\]
\end{enumerate}
\end{lemma}

\begin{proof}
(i) Evenness follows from the symmetry of $g$. The function $b\mapsto \mathbb{P}(|g+b|<a)$ is the Gaussian mass of an interval of length $2a$ centered at $-b$. By symmetry and strict unimodality of $\varphi$, it is strictly decreasing in $|b|$.

(ii) If $|g|<a$ and $|b|\ge 2a$, then $|g+b|\ge |b|-|g|>a$.

(iii) We have $G_a(b)=\Phi(-a-b)+1-\Phi(a-b)-\mathbb{P}(|g|\ge a)$, where $\Phi$ is the Gaussian distribution function. A direct computation gives $G_a'(0)=0$ and $G_a''(0)=2a\varphi(a)>0$. Hence the function $(a,b)\mapsto G_a(b)/b^2$, extended by $a\varphi(a)$ at $b=0$, is continuous and strictly positive on the compact set $[a_-,a_+]\times[0,1]$, by (i). So it is bounded below by a positive constant there. For $|b|\ge 1$, by (i) we have $G_a(b)\ge G_a(1)\ge\min_{a'\in[a_-,a_+]}G_{a'}(1)>0$. Since $a_\pm$ depend only on $c_1$, the constant $\gamma_0$ is absolute.
\end{proof}

The following lemma shows that a point with $m$ large coefficients is visible. This is the mechanism used in Section~\ref{sec:bayesian}.

\begin{lemma}[Strong coordinates]\label{lem:strong}
Let $t\in S^{n-1}$, $m\in\mathcal{M}$ and $r,u>0$ with $ru\ge 2a_m$. If
\[
\#\big\{j\in[N]:\ |\langle t,v_j\rangle|\ge u\big\}\ge m,
\]
then $t$ is visible at radius $r$.
\end{lemma}

\begin{proof}
Let $J:=\{j:|\langle t,v_j\rangle|\ge u\}$, so that $|J|\ge m$. For $j\in J$ we have $r|\langle t,v_j\rangle|\ge 2a_m$, hence $q_j\ge 1-p_m$ by Lemma~\ref{lem:gain}(ii). Therefore, using $p_m\le c_1\le 1/4$,
\[
\Delta_{m,J}\ge |J|(1-2p_m)\ge \frac m2.
\]
Moreover, $q_j(1-q_j)\le 1-q_j\le p_m$ for $j\in J$. Since $|J|\le N$ and $m^2\le N$,
\[
\sigma^2_{m,J}\le Np_m+|J|p_m\le 2Np_m=2c_1\frac{m^4}{N}\le 2c_1m^2.
\]
Thus $C_1\sigma_{m,J}\le C_1\sqrt{2c_1}\,m\le m/2\le\Delta_{m,J}$, because $c_1\le 1/(8C_1^2)$.
\end{proof}

The next corollary translates Lemma~\ref{lem:strong} into the logarithmic scale used in Sections~\ref{sec:lower-bounds} and~\ref{sec:bayesian}.

\begin{corollary}\label{cor:log-scale}
Let $t\in S^{n-1}$, let $\ell$ be an integer with $1\le \ell\le\lfloor\sqrt N\rfloor$, and let $r,u>0$. Suppose that
\[
\#\big\{j\in[N]:\ |\langle t,v_j\rangle|\ge u\big\}\ge \ell
\qquad\text{and}\qquad
r^2u^2\ge 4K\log\frac{eN}{\ell^2}.
\]
Then $t$ is visible at radius $r$.
\end{corollary}

\begin{proof}
Let $m\in\mathcal{M}$ satisfy $\ell/2<m\le\ell$. By \eqref{eq:threshold-size} and since $\log(eN/\ell^2)\ge 1$,
\[
a_m^2\le 4\log\frac{4N}{\ell^2}+2\log\frac1{c_1}\le \Big(4+4\log4+2\log\frac1{c_1}\Big)\log\frac{eN}{\ell^2}=K\log\frac{eN}{\ell^2}.
\]
Hence $ru\ge 2a_m$, and at least $\ell\ge m$ coefficients of $t$ have absolute value at least $u$. Lemma~\ref{lem:strong} applies.
\end{proof}

We can now make precise the relation between the Bayesian regularity of Definition~\ref{def:regular_bayes} and visibility. Let $C_{3.1}$ denote the absolute constant in Definition~\ref{def:regular_bayes}.

\begin{proposition}[Bayesian regularity means typical points are visible]\label{prop:bayes-visible}
Let $\mu$ be a probability measure supported on $T$ that is $(m,\delta)$-regular in the sense of Definition~\ref{def:regular_bayes}, with respect to the orthonormal family $V$. If $C_r^2\ge 12KC_{3.1}$, then
\[
\mu\Big(\Big\{t:\ t \text{ is visible at radius } \frac{C_r}{\sqrt{\kappa_\mu(T)}}\Big\}\Big)\ \ge\ 1-\delta.
\]
\end{proposition}

\begin{proof}
By Definition~\ref{def:regular_bayes}, with $\mu$-probability at least
$1-\delta$, the point $t$ has at least $m$ coefficients
of absolute value at least $u$, where $1\le m\le N$ and
\[
\kappa_\mu(T)\le
C_{3.1}\frac{u^2}{1+\log_+(N/m^2)}.
\]
Set $\ell=\min\{m,\lfloor\sqrt N\rfloor\}$.
Then $t$ has at least $\ell$ coefficients of
absolute value at least $u$. For $r=C_r/\sqrt{\kappa_\mu(T)}$, we obtain
\[
r^2u^2
\ge \frac{C_r^2}{C_{3.1}}
\left(1+\log_+\frac{N}{m^2}\right).
\]
Since
\[
\log\frac{eN}{\ell^2}
\le C_0\left(1+\log_+\frac{N}{m^2}\right)
\]
for some absolute constant $C_0$, choosing
$C_r^2\ge 4KC_0C_{3.1}$ gives
\[
r^2u^2\ge 4K\log\frac{eN}{\ell^2}.
\]
Choosing $C_0=3$, Corollary~\ref{cor:log-scale} concludes the proof.
\end{proof}

Proposition~\ref{prop:bayes-visible} and Definition~\ref{def:regular} make the comparison between the two problems transparent:
\begin{center}
\begin{tabular}{ll}
Bayesian problem: & $\mu$-most points of $T$ are visible at radius $\asymp \kappa_\mu(T)^{-1/2}$;\\[0.5ex]
Worst-Case problem: & \emph{every} point of $T$ is visible at radius $\asymp \kappa(T)^{-1/2}$.
\end{tabular}
\end{center}
In other words, the adversary gains nothing over a least favorable prior as long as it cannot find a point of $T$ that is harder for the count test to see than a typical point of that prior. 

\begin{corollary}\label{cor:bayes-equals-nonadv}
Let $T\subseteq S^{n-1}$ be $C_r$-regular, and let $\mu\in\mathcal{P}(T)$ satisfy $\kappa_\mu(T)\le C_0\,\kappa(T)$ for some $C_0\ge 1$. Then
\[
r_{\mu}(T)\ \asymp\ r(T)\ \asymp\ \frac{1}{\sqrt{\kappa(T)}},
\]
where the implied constants depend only on $C_0$ and $C_r$. In particular, since such $\mu$ exists with $C_0=2$,
\[
r(T)\ \asymp\ \sup_{\mu\in\mathcal{P}(T)} r_{\mu}(T).
\]
\end{corollary}

\begin{proof}
By Proposition~\ref{prop:lower_bound_bayesian} and Theorem~\ref{thm:regular},
\[
\frac{c}{\sqrt{C_0\kappa(T)}}\le\frac{c}{\sqrt{\kappa_\mu(T)}}\le r_{\mu}(T)\le r(T)\le\frac{C_r}{\sqrt{\kappa(T)}}. \qedhere
\]
\end{proof}

\begin{remark}
Suppose that for every $t\in T$ there exist $u>0$ and $1\le m\le\sqrt N$ with $\#\{j:|\langle t,v_j\rangle|\ge u\}\ge m$ and $\kappa(T)\le c_0u^2/\log(eN/m^2)$. Then, by Corollary~\ref{cor:log-scale}, $T$ is $C_r$-regular for every $C_r^2\ge 4Kc_0$. Thus Definition~\ref{def:regular} contains the strong-coordinate condition as a special case. The strong-coordinate condition alone, however, is not sufficient for $T_s$ when $s\gg\sqrt n$. For instance, a unit vector whose squared coefficients decay like $1/(j\log n)$ for $\sqrt n\le j\le n$ satisfies the strong-coordinate condition only up to a logarithmic term, while it is visible in the sense of Definition~\ref{def:visible} through the accumulation of many ``weak coordinates".
\end{remark}

\subsection{Sparse vectors}

Let $T_s=\{t\in S^{n-1}:|\operatorname{supp}(t)|\le s\}$, and take $V$ to be the canonical basis, so that $N=n$. By Section~\ref{s: sparse}, there is an absolute constant $C_\kappa\ge 1$ such that
\begin{equation}\label{eq:kappa-sparse-upper}
\kappa(T_s)\le\frac{C_\kappa}{s\log(en/s^2)}\quad\text{if } s\le\sqrt n,
\qquad
\kappa(T_s)\le\frac{C_\kappa}{\sqrt n}\quad\text{if } s>\sqrt n.
\end{equation}

We need the following elementary estimate on order statistics. It extends the claim in the proof of the regularity of $T_s$ in the sparse regime.

\begin{lemma}\label{lem:order-stats}
Let $x\in\mathbb{R}^n$ be supported on at most $k$ coordinates, where $k$ is an integer satisfying $1\le k\le\lfloor\sqrt n\rfloor$, and let $\eta:=\|x\|_2^2$. Then there exists $1\le\ell\le k$ such that
\[
|x|_{(\ell)}^2\ \ge\ \frac{\eta}{3}\cdot\frac{\log(en/\ell^2)}{k\log(en/k^2)}.
\]
\end{lemma}

\begin{proof}
If $\eta=0$, the conclusion is immediate. Assume $\eta>0$. Suppose, by contradiction, that $|x|_{(\ell)}^2<\theta\log(en/\ell^2)$ for every $1\le\ell\le k$, where $\theta:=\eta/(3k\log(en/k^2))$. Using $\log(k!)\ge k\log k-k$, we have
\[
\sum_{\ell=1}^{k}\log\frac{en}{\ell^2}=k\log(en)-2\log(k!)\le k\log\frac{en}{k^2}+2k\le 3k\log\frac{en}{k^2},
\]
where the last step uses $\log(en/k^2)\ge1$. Summing the assumed bound over $\ell$ yields $\eta=\sum_{\ell\le k}|x|_{(\ell)}^2<3\theta k\log(en/k^2)=\eta$, a contradiction.
\end{proof}

\begin{theorem}\label{thm:sparse-regular}
There exists an absolute constant $C_r>0$ such that $T_s$ is $C_r$-regular with respect to the canonical basis, for every $1\le s\le n$. Consequently,
\[
r(T_s)\ \asymp\ \frac{1}{\sqrt{\kappa(T_s)}}\ \asymp\
\begin{cases}
\sqrt{s\log(en/s^2)}, & s\le\sqrt n,\\[0.5ex]
n^{1/4}, & s>\sqrt n.
\end{cases}
\]
\end{theorem}

\begin{proof}
Fix $t\in T_s$ and let $r:=r_{T_s}=C_r/\sqrt{\kappa(T_s)}$, where $C_r$ will be chosen large enough, depending only on the absolute constants $C_1,c_1,K,\gamma_0,a_+$ and $C_\kappa$. We write $t_j=\langle t,e_j\rangle$.

\medskip
\noindent\emph{Sparse regime: $s\le\sqrt n$.}
Apply Lemma~\ref{lem:order-stats} with $x=t$, $k=s$ and $\eta=1$. This gives $1\le\ell\le s\le\sqrt n$ such that at least $\ell$ coordinates of $t$ have absolute value at least $u$, where
\[
u^2:=\frac{\log(en/\ell^2)}{3s\log(en/s^2)}.
\]
By \eqref{eq:kappa-sparse-upper}, $r^2\ge (C_r^2/C_\kappa)\,s\log(en/s^2)$. Hence
\[
r^2u^2\ge\frac{C_r^2}{3C_\kappa}\log\frac{en}{\ell^2}\ge 4K\log\frac{en}{\ell^2},
\]
provided $C_r^2\ge 12KC_\kappa$. By Corollary~\ref{cor:log-scale}, $t$ is visible at radius $r$.

\medskip
\noindent\emph{Dense regime: $s>\sqrt n$.}
By \eqref{eq:kappa-sparse-upper}, $r^2\ge A\sqrt n$, where $A:=C_r^2/C_\kappa$. Write $a^*:=a_{m^*}$, $p^*:=p_{m^*}$ for the top scale, and recall from \eqref{eq:top-scale} that $a^*\in[a_-,a_+]$ and $p^*\le c_1$. Split the coordinates into
\[
J_L:=\{j:\ r|t_j|\ge 2a^*\},\qquad J_S:=[n]\setminus J_L.
\]
Since $\sum_j t_j^2=1$, one of the following three cases holds.

\smallskip
\noindent\emph{Case 1: $\sum_{j\in J_S}t_j^2\ge 1/2$.} We take $m=m^*$ and $J=J_S$. For $j\in J_S$ we have $r|t_j|<2a^*\le 2a_+$, so that $\min\{r^2t_j^2,1\}\ge\beta_0r^2t_j^2$ with $\beta_0:=1/\max\{1,4a_+^2\}$. By Lemma~\ref{lem:gain}(iii),
\[
\Delta:=\Delta_{m^*,J_S}\ \ge\ \gamma_0\beta_0r^2\sum_{j\in J_S}t_j^2\ \ge\ \frac{\gamma_0\beta_0A}{2}\sqrt n.
\]
Moreover, $q_j(1-q_j)\le q_j=p^*+(q_j-p^*)$, so
\[
\sigma^2_{m^*,J_S}\le Np^*+|J_S|p^*+\Delta\le 2c_1n+\Delta.
\]
If $A$ is large enough that $\gamma_0\beta_0A\ge \max\{1,\,4C_1\sqrt{c_1},\,4C_1^2\}$, then $\Delta\ge\sqrt n/2\ge m^*/2$, $\Delta^2\ge 4C_1^2c_1n$ and $\Delta\ge 2C_1^2$. The last two inequalities give $C_1^2(2c_1n+\Delta)\le\Delta^2/2+\Delta^2/2$, that is, $\Delta\ge C_1\sigma_{m^*,J_S}$. Hence \eqref{eq:visible} holds.

\smallskip
\noindent\emph{Case 2: $\sum_{j\in J_L}t_j^2\ge 1/2$ and $|J_L|\ge m^*$.} At least $m^*$ coordinates satisfy $|t_j|\ge u:=2a^*/r$, and $ru=2a_{m^*}$. Lemma~\ref{lem:strong} shows that $t$ is visible at radius $r$.

\smallskip
\noindent\emph{Case 3: $\sum_{j\in J_L}t_j^2\ge 1/2$ and $k:=|J_L|<m^*\le\sqrt n$.} Let $x$ be the restriction of $t$ to $J_L$. Then $x$ is supported on $k\le\sqrt n$ coordinates and $\|x\|_2^2\ge1/2$. For $1\le k\le\sqrt n$ we have
\[
k\log\frac{en}{k^2}=\sqrt n\,\frac{1+2\log y}{y}\Big|_{y=\sqrt n/k}\le\frac{2}{\sqrt e}\sqrt n\le 2\sqrt n.
\]
So Lemma~\ref{lem:order-stats} yields $1\le\ell\le k$ such that at least $\ell$ coordinates of $t$ have absolute value at least $u$, where
\[
u^2:=\frac{\log(en/\ell^2)}{12\sqrt n}.
\]
Then $r^2u^2\ge (A/12)\log(en/\ell^2)\ge 4K\log(en/\ell^2)$ if $A\ge 48K$. Since $\ell\le\sqrt n$, Corollary~\ref{cor:log-scale} shows that $t$ is visible at radius $r$.

\medskip
Choosing $C_r$ so large that all the above requirements on $C_r$ and $A=C_r^2/C_\kappa$ hold, every $t\in T_s$ is visible at radius $r_{T_s}$, for every $1\le s\le n$. The rates follow from Theorem~\ref{thm:regular} and the estimates on $\kappa(T_s)$ of Section~\ref{s: sparse}.
\end{proof}

\begin{remark}[Open problem]
Regularity is a sufficient condition. A complete characterization of the sets $T$ for which the lower bound $r(T)\gtrsim\kappa(T)^{-1/2}$ is sharp remains open. In particular, it would be interesting to know if $\kappa(T)^{-1/2}$ is the right quantity to capture the order of the critical radius $r(T)$.
\end{remark}

\subsection{Another Application: Sobolev ellipsoids}\label{subsec:sobolev}

We now consider a class that is not sparse in any sense: signals whose energy decays with frequency. For $\beta>0$ and $1\le L\le n^{\beta}$, let
\[
T_{\beta,L}:=\Big\{t\in S^{n-1}:\ \sum_{j=1}^n j^{2\beta}t_j^2\le L^2\Big\},
\qquad
D:=\big\lfloor L^{1/\beta}\big\rfloor\in[1,n].
\]
Here the coordinates play the role of the coefficients of a function in a Fourier-type basis, and the constraint is the usual Sobolev ellipsoid of smoothness $\beta$ and radius $L$. Typical elements of $T_{\beta,L}$ are dense: they spread their energy over roughly $D$ coordinates with decaying amplitudes. Testing a Gaussian mean against smooth alternatives of this type is the classical problem of nonparametric goodness-of-fit testing, initiated by Ermakov and Ingster (see \cite{Ermakov1990,Ingster1993I} and the references therein).

The main point of this subsection is that both regularity conditions are easy to verify for $T_{\beta,L}$, and that the resulting rate recovers the classical minimax rate. The verification relies on the following general criterion, which isolates the argument used in the dense regime of Theorem~\ref{thm:sparse-regular}.

\begin{proposition}[Dense-regime criterion]\label{prop:dense-criterion}
Let $V=\{v_1,\dots,v_N\}$ be an orthonormal family, and suppose that
\begin{equation}\label{eq:energy-in-V}
\sum_{j=1}^N\langle t,v_j\rangle^2\ \ge\ \frac12\qquad\text{for every } t\in T.
\end{equation}
There is an absolute constant $A_0>0$ such that every $t\in T$ is visible at every radius $r$ with $r^2\ge A_0\sqrt N$. In particular,
\begin{enumerate}
\item[(i)] $r(T)\le A_0^{1/2}N^{1/4}$;
\item[(ii)] if $\kappa(T)\le C_0N^{-1/2}$, then $T$ is $C_r$-regular in the sense of Definition~\ref{def:regular}, with $C_r=(A_0C_0)^{1/2}$.
\end{enumerate}
\end{proposition}

\begin{proof}
Fix $t\in T$ and $r$ with $r^2\ge A\sqrt N$, where $A$ is a large absolute constant to be chosen. Write $c_j=\langle t,v_j\rangle$, and let $m^*$, $a^*:=a_{m^*}$ and $p^*:=p_{m^*}$ denote the top scale, so that $a^*\in[a_-,a_+]$ and $p^*\le c_1$ by \eqref{eq:top-scale}. Split $[N]$ into
\[
J_L:=\{j:\ r|c_j|\ge 2a^*\},\qquad J_S:=[N]\setminus J_L.
\]
By \eqref{eq:energy-in-V}, one of the following three cases occurs.

\smallskip
\noindent\emph{Case 1: $\sum_{j\in J_S}c_j^2\ge 1/4$.} Take $m=m^*$ and $J=J_S$. For $j\in J_S$ we have $r|c_j|<2a_+$, so $\min\{r^2c_j^2,1\}\ge\beta_0r^2c_j^2$ with $\beta_0:=1/\max\{1,4a_+^2\}$. By Lemma~\ref{lem:gain}(iii),
\[
\Delta:=\Delta_{m^*,J_S}\ \ge\ \gamma_0\beta_0r^2\sum_{j\in J_S}c_j^2\ \ge\ \frac{\gamma_0\beta_0A}{4}\sqrt N.
\]
Since $q_j(1-q_j)\le q_j=p^*+(q_j-p^*)$, we have $\sigma^2_{m^*,J_S}\le 2Np^*+\Delta\le 2c_1N+\Delta$. If $\gamma_0\beta_0A\ge\max\{2,\,8C_1\sqrt{c_1},\,8C_1^2\}$, then $\Delta\ge\sqrt N/2\ge m^*/2$, $\Delta^2\ge 4C_1^2c_1N$ and $\Delta\ge 2C_1^2$. The last two inequalities give $C_1^2(2c_1N+\Delta)\le\Delta^2$, so \eqref{eq:visible} holds.

\smallskip
\noindent\emph{Case 2: $\sum_{j\in J_L}c_j^2\ge 1/4$ and $|J_L|\ge m^*$.} At least $m^*$ coefficients satisfy $|c_j|\ge u:=2a^*/r$, and $ru=2a_{m^*}$. Lemma~\ref{lem:strong} applies.

\smallskip
\noindent\emph{Case 3: $\sum_{j\in J_L}c_j^2\ge 1/4$ and $k:=|J_L|<m^*\le\sqrt N$.} Apply Lemma~\ref{lem:order-stats} in $\mathbb{R}^N$ to the vector $x=(c_j\mathbf 1_{j\in J_L})_{j\le N}$, which is supported on $k\le\sqrt N$ coordinates and satisfies $\|x\|_2^2\ge1/4$. As in the proof of Theorem~\ref{thm:sparse-regular}, $k\log(eN/k^2)\le 2\sqrt N$ for $1\le k\le\sqrt N$. We obtain $1\le\ell\le k$ such that at least $\ell$ coefficients satisfy $|c_j|\ge u$, where
\[
u^2:=\frac{\log(eN/\ell^2)}{24\sqrt N}.
\]
Then $r^2u^2\ge (A/24)\log(eN/\ell^2)\ge 4K\log(eN/\ell^2)$ if $A\ge 96K$. Since $\ell\le\sqrt N$, Corollary~\ref{cor:log-scale} applies.

\smallskip
Taking $A_0:=\max\{(\gamma_0\beta_0)^{-1}\max\{2,8C_1\sqrt{c_1},8C_1^2\},\,96K\}$, every $t\in T$ is visible at every radius $r$ with $r^2\ge A_0\sqrt N$. Part (i) then follows from the proof of Theorem~\ref{thm:regular}, which only uses visibility at the given radius. For (ii), $r_T^2=C_r^2/\kappa(T)\ge (C_r^2/C_0)\sqrt N=A_0\sqrt N$.
\end{proof}

We now verify the two regularity conditions for $T_{\beta,L}$. Throughout, constants depending only on $\beta$ are denoted by $C_\beta,c_\beta$.

\begin{lemma}[Energy at low frequencies]\label{lem:sobolev-tail}
Let $N:=\min\{n,\lceil 2^{1+1/(2\beta)}D\rceil\}$. Then every $t\in T_{\beta,L}$ satisfies
\[
\sum_{j\le N}t_j^2\ \ge\ \frac12.
\]
\end{lemma}

\begin{proof}
If $N=n$ there is nothing to prove. Otherwise, since $L^{1/\beta}<D+1\le 2D$, we have $N^{2\beta}\ge 2^{2\beta+1}D^{2\beta}\ge 2L^2$. Hence
\[
\sum_{j>N}t_j^2\ \le\ N^{-2\beta}\sum_{j>N}j^{2\beta}t_j^2\ \le\ \frac{L^2}{N^{2\beta}}\ \le\ \frac12. \qedhere
\]
\end{proof}

\begin{proposition}[Bayesian regularity]\label{prop:sobolev-bayes}
Let $c>0$ be the absolute constant of the dense regime in Section~\ref{sec:bayesian}, and let $m:=\max\{1,\lfloor c\sqrt D\rfloor\}$. Let $\mu$ be the uniform measure on the flat $m$-sparse vectors supported on the first $D$ coordinates,
\[
\Big\{\frac{1}{\sqrt m}\sum_{i\in I}\varepsilon_ie_i:\ I\subseteq[D],\ |I|=m,\ \varepsilon_i\in\{\pm1\}\Big\}.
\]
Then $\mu$ is supported on $T_{\beta,L}$, it is $(m,0)$-regular in the sense of Definition~\ref{def:regular_bayes} with respect to $V=\{e_1,\dots,e_D\}$, and
\[
\kappa_\mu(T_{\beta,L})\ \asymp\ D^{-1/2}.
\]
\end{proposition}

\begin{proof}
For $t$ in the support of $\mu$,
\[
\sum_{j}j^{2\beta}t_j^2\ \le\ D^{2\beta}\sum_{j\le D}t_j^2\ =\ D^{2\beta}\ \le\ L^2,
\]
so $t\in T_{\beta,L}$. Every such $t$ has exactly $m$ coefficients of absolute value $u:=m^{-1/2}$ with respect to $V$, and $m\le\sqrt D$. The measure $\mu$ is the flat sparse prior of Section~\ref{sec:bayesian} in dimension $D$, and the computation of the dense regime there (with $n$ replaced by $D$) gives
\[
\kappa_\mu\ \asymp\ \frac{1}{m\log(eD/m^2)}\ \asymp\ D^{-1/2}\ \asymp\ \frac{u^2}{\log(eD/m^2)}.
\]
When $c\sqrt D<1$, so that $m=1$ and $D$ is bounded, all quantities are of order one and the same conclusion holds. Thus both conditions of Definition~\ref{def:regular_bayes} hold with $\delta=0$.
\end{proof}

\begin{theorem}\label{thm:sobolev}
Let $\beta>0$ and $1\le L\le n^\beta$. Then $T_{\beta,L}$ satisfies the regularity of Section~\ref{sec:bayesian} (Proposition~\ref{prop:sobolev-bayes}), and it is $C_\beta$-regular in the sense of Definition~\ref{def:regular} with respect to $V=\{e_1,\dots,e_N\}$, where $N$ is as in Lemma~\ref{lem:sobolev-tail}. Moreover,
\[
r(T_{\beta,L})^2\ \asymp_\beta\ r_{\mu}(T_{\beta,L})^2\ \asymp_\beta\ \frac{1}{\kappa(T_{\beta,L})}\ \asymp_\beta\ \sqrt D\ \asymp_\beta\ L^{1/(2\beta)}.
\]
\end{theorem}

\begin{proof}
By Lemma~\ref{lem:sobolev-tail}, condition \eqref{eq:energy-in-V} holds with $V=\{e_1,\dots,e_N\}$. Moreover $N\le C_\beta D$, so Proposition~\ref{prop:sobolev-bayes} gives
\[
\kappa(T_{\beta,L})\ \le\ \kappa_\mu\ \le\ C\,D^{-1/2}\ \le\ C_\beta N^{-1/2}.
\]
Proposition~\ref{prop:dense-criterion}(ii) shows that $T_{\beta,L}$ is $C_\beta$-regular. Next, Proposition~\ref{prop:dense-criterion}(i) gives $r(T_{\beta,L})\lesssim N^{1/4}\lesssim_\beta D^{1/4}$. Conversely, by Proposition \ref{prop:lower_bound_bayesian} and Proposition~\ref{prop:sobolev-bayes},
\[
r(T_{\beta,L})\ \ge\ r_{\mu}(T_{\beta,L})\ \gtrsim\ \kappa_\mu^{-1/2}\ \gtrsim\ D^{1/4}.
\]
Finally, Proposition~\ref{prop:lower_bound_bayesian} and the upper bound give $\kappa(T_{\beta,L})^{-1/2}\lesssim r(T_{\beta,L})\lesssim_\beta D^{1/4}$, while $\kappa(T_{\beta,L})\le\kappa_\mu\lesssim D^{-1/2}$. Hence $\kappa(T_{\beta,L})\asymp_\beta D^{-1/2}$, and $D\asymp L^{1/\beta}$.
\end{proof}

\section{Declaration of AI Use}
The authors used AI for proofreading and computations. The main ideas in this project were developed prior the emergence of proofs generated by the recent generative models.

\section{Acknowledgments}
R.V is supported by NSF Grant DMS 2451011 and U.S. Air Force Grant FA9550-25-1-0294.

\end{document}